\documentclass[epsfig,latexsym,amsfonts,twoside]{article}
\usepackage{amssymb,amsmath,amscd,epsfig}

\def\part#1{\frac{\partial\phantom{#1}}{\partial#1}}
\newtheorem{thm}{Theorem}

\newtheorem{lem}[thm]{Lemma}
\newtheorem{cor}[thm]{Corollary}
\newtheorem{cnj}[thm]{Conjecture}
\newenvironment{prf}{\begin{trivlist}\item[]{\bf Proof} }%
{\hfill $\Box$ \end{trivlist}}
\newenvironment{dfn}{\begin{trivlist}\item[]{\bf Definition}\em }%
{\end{trivlist}}
\newenvironment{rmk}{\begin{trivlist}\item[]{\bf Remark} }%
{\end{trivlist}}
\newenvironment{exm}{\begin{trivlist}\item[]{\bf Example} }%
{\end{trivlist}}

\def\Z{\ifmmode{{\mathbb Z}}\else{${\mathbb Z}$}\fi}
\def\Q{\ifmmode{{\mathbb Q}}\else{${\mathbb Q}$}\fi}
\def\C{\ifmmode{{\mathbb C}}\else{${\mathbb C}$}\fi}
\def\P{\ifmmode{{\mathbb P}}\else{${\mathbb P}$}\fi}

\def\H{\ifmmode{{\mathrm H}}\else{${\mathrm H}$}\fi}

\def\B{\ifmmode{{\cal B}}\else{${\cal B}$}\fi}
\def\E{\ifmmode{{\cal E}}\else{${\cal E}$}\fi}
\def\F{\ifmmode{{\cal F}}\else{${\cal F}$}\fi}
\def\K{\ifmmode{{\cal K}}\else{${\cal K}$}\fi}
\def\L{\ifmmode{{\cal L}}\else{${\cal L}$}\fi}
\def\M{\ifmmode{{\cal M}}\else{${\cal M}$}\fi}
\def\N{\ifmmode{{\cal N}}\else{${\cal N}$}\fi}
\def\O{\ifmmode{{\cal O}}\else{${\cal O}$}\fi}
\def\U{\ifmmode{{\cal U}}\else{${\cal U}$}\fi}
\def\X{\ifmmode{{\cal X}}\else{${\cal X}$}\fi}

\def\Br{\ifmmode{{\mathrm{Br}}}\else{${\mathrm{Br}}$}\fi}
\def\OG{\ifmmode{\widetilde{\cal M}_4}\else{$\widetilde{\cal M}_4$}\fi}
\def\D{\ifmmode{{\cal{D}}_{\mathrm{coh}}^b}\else{${{\cal{D}}_{\mathrm{coh}}^b}$}\fi}
\def\Shah{\ifmmode{\amalg\hspace*{-3.5pt}\amalg}\else{$\amalg\hspace*{-3.5pt}\amalg$}\fi}

\begin{document}

\title{On Lagrangian fibrations by Jacobians I\footnote{2000 {\em Mathematics Subject
    Classification.\/} 14J28; 14D06, 53C26.}}
\author{Justin Sawon}
\date{September, 2011}
\maketitle

\begin{abstract}
Let $Y\rightarrow\P^n$ be a flat family of integral Gorenstein curves,
such that the compactified relative Jacobian $X=\overline{J}^d(Y/\P^n)$ is a Lagrangian fibration. We prove that the degree of the discriminant locus $\Delta\subset\P^n$ is at least $4n+2$, and we prove that $X$ is a Beauville-Mukai integrable system if $\mathrm{deg}\Delta>4n+20$.
\end{abstract}

\section{Introduction}

An irreducible holomorphic symplectic manifold is a compact,
simply-connected, K{\"a}hler manifold $X$ which admits a
non-degenerate holomorphic two-form $\sigma$ which generates
$\H^0(X,\Omega^2)$. Non-degeneracy means that $\sigma$ induces an
isomorphism $T\cong \Omega^1$, or equivalently, $\sigma^{\wedge n}$
trivializes the canonical bundle $K_X$, where $\mathrm{dim}X=2n$. Let
$\pi:X\rightarrow B$ be a proper surjective morphism with connected
fibres, with $B$ smooth and
$0<\mathrm{dim}B<2n$. Matsushita~\cite{matsushita99, matsushita00,
  matsushita05} proved that the base and fibres must have dimension
$n$, the fibres must be Lagrangian with respect to the holomorphic
symplectic form $\sigma$, and the generic fibre must be an abelian variety. His
proof assumes that $X$ is projective, but can be easily adapted to the
non-projective case (see Huybrechts, Proposition 24.8~\cite{ghj02}). In the projective case, Hwang~\cite{hwang08} proved that
the base $B$ must be isomorphic to $\P^n$. We expect $B\cong\P^n$ to
be true also in the non-projective case, and we will call such a
fibration $\pi:X\rightarrow\P^n$ a {\em Lagrangian fibration\/}.

In~\cite{sawon03} the author described how one might use Lagrangian
fibrations to classify holomorphic symplectic manifolds up to
deformation. All of the known examples of holomorphic symplectic
manifolds can be deformed to Lagrangian fibrations (see
Beauville~\cite{beauville99}, Debarre~\cite{debarre99}, and
Rapagnetta~\cite{rapagnetta07}). In fact, it is expected that
Lagrangian fibrations will be dense in the moduli space of all
deformations, and this has been proved in some cases (see
Gulbrandsen~\cite{gulbrandsen07}, Markushevich~\cite{markushevich06},
Yoshioka~\cite{yoshioka09}, and the author's article~\cite{sawon07}).

In this article we work towards a classification of Lagrangian fibrations by Jacobians. Examples of such fibrations are given by the Beauville-Mukai integrable systems~\cite{beauville99, mukai84}: one starts with a family of genus $n$ curves $Y/\P^n$ which is a complete linear system of curves in a K3 surface $S$, and constructs its compactified relative Jacobian $X=\overline{J}^d(Y/\P^n)$. Then $X\rightarrow\P^n$ is a Lagrangian fibration and $X$ is a deformation of $\mathrm{Hilb}^nS$ (see Section~2.3 for more details).

Conversely, we expect the following to be true.
\begin{cnj}
\label{conjecture}
Let $X\rightarrow\P^n$ be a Lagrangian fibration whose fibres are Jacobians of genus $n$ curves. Then $X$ is a Beauville-Mukai integrable system.
\end{cnj}
We first need to clarify what it means to be a Lagrangian fibration by Jacobians. We start with a family $Y/\P^n$ of geometrically integral (reduced and irreducible) Gorenstein curves. The compactified Jacobian is well-defined for integral curves (D'Souza~\cite{dsouza79}, Altman and Kleiman~\cite{ak80}). We let
$$X:=\overline{J}^d(Y/\P^n)$$
and our hypothesis is that $X$ is a Lagrangian fibration over $\P^n$. We prove (Theorem~\ref{main})
\begin{enumerate}
\item the degree of the discriminant locus $\Delta\subset\P^n$, which parametrizes singular fibres, is at least $4n+2$,
\item if $\mathrm{deg}\Delta>4n+20$ then $X$ is a Beauville-Mukai integrable system, i.e., the family of curves $Y/\P^n$ is a complete linear system of curves on a K3 surface.
\end{enumerate}
We also verify the bound on $\mathrm{deg}\Delta$ when $n=2$, thereby giving a new proof of a result (Corollary~\ref{corollary}) originally due to Markushevich~\cite{markushevich96}.

Our proof uses a construction of Hurtubise~\cite{hurtubise96}. We first show that it suffices to prove the theorem for degree $d=1$ Jacobians, in which case the Abel map gives a canonical embedding $Y\hookrightarrow X$. If the restriction $\sigma|_Y$ of the holomorphic symplectic form $\sigma$ to $Y$ has rank two, then the null directions of $\sigma|_Y$ lead to a rank $n-1$ foliation $F$ known as the {\em characteristic foliation\/}. Moreover, the
two-form $\sigma|_Y$ descends to a non-degenerate two-form on the space of leaves
$Q=Y/F$. In this way one obtains a holomorphic symplectic surface $Q$
from the Lagrangian fibration. The curves in the family $Y/\P^n$ project down to $Q$.

Hurtubise's argument is local: he starts with a family of smooth curves
parametrized by a ball in $\C^n$. The complex surface $Q$ that he obtains is non-compact, though for genus $n\geq 3$ it admits a projective embedding as an open subset (in the analytic topology) of an algebraic surface. The principal
difficulty in applying his argument in a global setting is in dealing
with the characteristic foliation. If we want the space of leaves to be reasonably
behaved (for example, Hausdorff) we must first show that the foliation
is algebraically integrable, which in particular implies that it has
compact leaves. The earliest results in this direction were obtained by
Miyaoka~\cite{miyaoka87}, and later developed by Bogomolov and
McQuillan; we will
follow more closely a recent refinement by Kebekus, Sol{\'a} Conde,
and Toma~\cite{ksct07}. The key idea is that if $Y$ is covered by
curves on which $F$ is ample, then by applying Mori's bend-and-break
argument one can produce rational curves which must be contained in
the leaves of $F$. One thereby shows that the leaves are rationally
connected, and in particular algebraic. In our case, we are able to find the required rational curves in $Y$ more directly.

Hurtubise also has to assume that $\sigma|_Y$ has rank two, which implies
that $Y$ is a coisotropic submanifold of $X$. In our case, we prove that the 
rank of $\sigma|_Y$ must equal two. In real symplectic geometry the idea of 
taking quotients of coisotropic submanifolds is not new; however, this 
construction has rarely been used in holomorphic symplectic geometry (although, see Hwang and Oguiso~\cite{ho09}). A large number of examples of holomorphic coisotropic reduction, and some potential applications, are described in the author's article~\cite{sawon09}.

In a sequel~\cite{sawon11} to this paper we will prove Conjecture~\ref{conjecture} in a number of low dimensional cases: when $n=3$, $4$, and $5$ and all of the curves in the family $Y/\P^n$ are non-hyperelliptic, and when $n=3$ and all of the curves are hyperelliptic.

The author would like to thank Fabrizio Catanese, Brendan Hassett,
Jun-Muk Hwang, Stefan Kebekus, Manfred Lehn, Dimitri Markushevich,
Rick Miranda, Keiji Oguiso, and Christian Thier for many helpful
discussions on the material presented here, and is grateful for the
hospitality of the Max-Planck-Institut f{\"u}r Mathematik, Bonn, and
the Institute for Mathematical Sciences, the Chinese University of
Hong Kong, where these results were obtained.

\section{Preliminaries}

\subsection{Compactified Jacobians}

Let $C$ be a geometrically integral (reduced and irreducible) curve of arithmetic genus $n$, and let $\tilde{C}$ be the normalization of $C$. The Jacobian $J$ of $C$ is the group scheme parametrizing rank-one locally-free sheaves on $C$. Here and throughout, we use $J$ to denote a single connected component of the Jacobian, omitting reference to the degree $d$ unless necessary; the same convention will apply to $\tilde{J}$, $\overline{J}$, etc. The Jacobian $J$ is an iterated extension of the Jacobian $\tilde{J}$ of $\tilde{C}$ by copies of $\C^*$ and $\C$. There is a natural compactification of $J$ given by the moduli space $\overline{J}$ of rank-one torsion-free sheaves on $C$. For integral curves this moduli space was constructed by D'Souza~\cite{dsouza79} (see also Altman and Kleiman~\cite{ak80}). If $\overline{J}$ is irreducible, its normalization $P$ will be an iterated $\P^1$-bundle over the abelian variety $\tilde{J}$, with the $\P^1$ fibres arising as compactifications of the $\C^*$ and $\C$ fibres; $P$ is the scheme representing the Presentation Functor (see Altman and Kleiman~\cite{ak90}). In general, $\overline{J}$ need not be irreducible (see Altman, Iarrobino, and Kleiman~\cite{aik77}), but a similar statement can be made about the closure of $J$ in $\overline{J}$ instead.

\begin{exm}
The compactified Jacobian of a curve with a single simple node was described by Oda and Seshadri~\cite{os79}. Suppose that the node of $C$ is obtained by identifying points $p$ and $q$ in $\tilde{C}$. Then $P$ is the $\P^1$-bundle $\P(\mathcal{L}_p\oplus\mathcal{L}_q)$ over $\tilde{J}$, where $\mathcal{L}_p$ and $\mathcal{L}_q$ are translates of the Poincar{\'e} bundle over $\tilde{J}$. Note that $P$ has distinguished sections $s_0$ and $s_{\infty}$. The compactified Jacobian $\overline{J}$ is obtained by gluing $s_0$ to $s_{\infty}$, but we include a translation by $\O(p-q)\in\tilde{J}$ when we identify $s_0\cong\tilde{J}$ and $s_{\infty}\cong\tilde{J}$. 
\end{exm}

\begin{exm}
The compactified Jacobian of a curve with an ordinary cusp was described by Kleiman~\cite{kleiman84} (see also Altman and Kleiman~\cite{ak90}). Suppose that the preimage in $\tilde{C}$ of the cusp is the point $p$. Then $\pi:P\rightarrow\tilde{J}$ is the $\P^1$-bundle $\P(J^1\mathcal{L}_p)$, where $J^1\mathcal{L}_p$ is the first jet bundle of $\mathcal{L}_p$. This time $P$ has a single distinguished section $s_{\infty}$ (the $\P^1$ fibres are obtained by compactifying $\C$, rather than $\C^*$ as in the previous example). There is also a particular vector field $v$ along $s_{\infty}$, i.e., a section of $TP|_{s_{\infty}}$. This vector field points out of $s_{\infty}$, and if we project down to $\tilde{J}$, we get a vector field $d\pi(v)$ which at each point $q\in\tilde{J}$ points in the direction of the curve $\tilde{C}$, if we take the Abel embedding of $\tilde{C}$ in $\tilde{J}$ and then translate it so that $p\in\tilde{C}$ is mapped to $q\in\tilde{J}$. The main point here is that $v$ is {\em not\/} tangent to the $\P^1$ fibres of $P$.

The compactified Jacobian $\overline{J}$ is the image of a birational morphism $f:P\rightarrow\overline{J}$ which contracts $v$. This creates a family of cusps along the image $\mathrm{Sing}\overline{J}$ of $s_{\infty}$: locally $\overline{J}$ looks like the product of a cusp and $\C^{n-1}$. Let $(z_1,z_2,\ldots,z_n)$ be local coordinates on $P$, with $s_{\infty}$ given by $z_1=0$ and $v$ given by $\frac{\partial\phantom{z}}{\partial z_1}+\frac{\partial\phantom{z}}{\partial z_2}$. Locally, the map to $\overline{J}$ is given by
$$f(z_1,z_2,\ldots,z_n)=(z_1^2,z_1^3,z_2-z_1,z_3,\ldots,z_n)\in \{(x,y,z,\ldots )\in\C^{n+1}|y^2=x^3\}.$$
Note that a $\P^1$ fibre of $P$, which is given by keeping $z_2,\ldots,z_n$ constant, will map isomorphically to its image in $\overline{J}$, i.e., the image is {\em not\/} a cuspidal curve, but a smooth rational curve which is tangent to $\mathrm{Sing}\overline{J}$ at the point where they meet.
\end{exm}

\subsection{Lagrangian fibrations}

Let $Y\rightarrow B$ be a flat family of geometrically integral (reduced and irreducible) curves over a projective manifold $B$. Let
$$X=\overline{J}^d(Y/B)$$
be the compactified relative Jacobian of $Y/B$~\cite{ak80}. Our main hypothesis throughout will be
$$\mbox{\em{suppose that }}\pi:X\rightarrow B\mbox{\em{ is a Lagrangian fibration.}}$$
In other words, $X$ is an irreducible holomorphic symplectic manifold of dimension $2n$ (we assume $2n\geq 4$, the two-dimensional case being well understood), and the fibres of the morphism $\pi$ are Lagrangian with respect to the holomorphic symplectic form.

\begin{rmk}
In particular, $X$ is smooth by hypothesis. In the $n=2$ case, Markushevich~\cite{markushevich96} found necessary and sufficient conditions on a family of genus two curves to ensure that the compactified relative Jacobian $X=\overline{J}^d(Y/B)$ be smooth. Since we are only concerned with Lagrangian fibrations in this paper, we will assume from the outset that $X$ is also a Lagrangian fibration, and use this to deduce some properties of the curves.
\end{rmk}

\begin{lem}
\label{genusn}
Every curve in the family $Y/B$ has arithmetic genus $n$.
\end{lem}

\begin{prf}
If $C$ is a curve in the family $Y/B$, then the dimension of its Jacobian $J$ is equal to the arithmetic genus of $C$. The closure of $J$ inside the compactified Jacobian $\overline{J}$ of $C$ forms an irreducible component of $\overline{J}$. Matsushita~\cite{matsushita00} proved that every irreducible component of a fibre of the fibration $\pi:X\rightarrow B$ must be Lagrangian, and in particular, must have dimension $n$.
\end{prf}

\begin{lem}
\label{base}
The base $B$ is isomorphic to $\P^n$.
\end{lem}

\begin{prf}
The base is smooth and projective by hypothesis. Therefore $Y$ and $X$ are also projective, since the curves have genus $n\geq 2$. Hwang~\cite{hwang08} proved that the base of a projective Lagrangian fibration with smooth and projective base must be isomorphic to $\P^n$.
\end{prf}

\begin{lem}
The generic curve in the family $Y/B$ is a smooth genus $n$ curve.
\end{lem}

\begin{prf}
By the holomorphic Liouville Theorem the generic fibre of $\pi:X\rightarrow B$ is a smooth $n$-dimensional abelian variety. It follows that the generic curve must be smooth, since it is integral.
\end{prf}

Irreducibility of the curves is essential for the above lemma, as the following example shows.

\begin{exm}
Let $S\rightarrow\P^1$ be an elliptic K3 surface. There is an induced morphism $\mathrm{Hilb}^nS\rightarrow\mathrm{Sym}^n\P^1\cong\P^n$ which makes $\mathrm{Hilb}^nS$ into a Lagrangian fibration. Over a generic point $\{p_1,\ldots,p_n\}\in\mathrm{Sym}^n\P^1$ of the base, the fibre of the Lagrangian fibration is $E_1\times\ldots\times E_n$, where $E_i$ is the elliptic fibre of $S\rightarrow\P^1$ above the point $p_i$. Of course $E_1\times\ldots\times E_n$ is smooth, but we can think of it as the Jacobian of the singular curve
$$C=E_1*E_2*\cdots *E_n$$
formed by joining $E_i$ to $E_{i+1}$ at a single point, for each $i=1,\ldots,n-1$.
\end{exm}

\begin{lem}
\label{singlenode}
There is a hypersurface $\Delta\subset B$ parametrizing singular fibres of $\pi:X\rightarrow B$ (equivalently, singular curves in the family $Y/B$). A curve above a generic point of $\Delta$ will contain a single simple node.
\end{lem}

\begin{prf}
Suppose $\Delta$ is empty. Then the homotopy long exact sequence of $X\rightarrow B$ gives
$$\ldots\rightarrow\pi_2(X)\rightarrow\pi_2(B)\rightarrow\pi_1(F)\rightarrow\pi_1(X)\rightarrow\pi_1(B)\rightarrow 1.$$
By Lemma~\ref{base}, $B\cong\P^n$, so $\pi_2(B)\cong\Z$ and $\pi_1(B)$ is trivial. The fibre $F$ is  a torus of real dimension $2n$, so $\pi_1(F)\cong\Z^{\oplus 2n}$. Then $\pi_1(X)$ would be non-trivial, contradicting the irreducibility of $X$. It follows that $\Delta$ is non-empty.

The fact that $\Delta$ is a hypersurface if it is non-empty is a non-trivial statement, proved by Hwang and Oguiso, Proposition 3.1(2) of~\cite{ho09}.

Next we use Hwang and Oguiso's classification of generic singular fibres of Lagrangian fibrations. Let $Z$ be an irreducible component of a generic singular fibre, with normalization $\tilde{Z}$. Then by Theorem 1.3 of~\cite{ho09}, the Albanese variety of $\tilde{Z}$ has dimension $n-1$. In the case of the compactified Jacobian of a curve $C$, we let $Z$ be the closure of $J$ in $\overline{J}$. Then the normalization is the scheme $P$ representing the Presentation Functor, and its Albanese variety is the Jacobian $\tilde{J}$ of the normalization $\tilde{C}$ of $C$. It follows that $\tilde{C}$ has genus $n-1$, $\delta=1$, and $C$ has either a single node or a single cusp.

Suppose that $C$ has a single cusp. In the terminology of~\cite{ho09}, the image of a $\P^1$ fibre on $P$ will be a {\em characteristic $1$-cycle\/} on $\overline{J}$. As we saw in the previous subsection, each $\P^1$ fibre maps isomorphically to its image in $\overline{J}$, which would mean that the characteristic $1$-cycles are smooth rational curves. This would contradict Theorem 1.4 of~\cite{ho09} (see also Propositions 4.4(3), 4.7(3), or 4.12(1)).
\end{prf}

\begin{rmk}
Let us give a heuristic explanation of why a cuspidal curve cannot appear in codimension one. Choose local coordinates $(z_1,z_2,\ldots,z_n)$ on $B$ such that $\Delta$ is given by $z_1=0$. The characteristic foliation on $\pi^{-1}(\Delta)$ is generated by the vector field $v_1$ dual (with respect to the holomorphic symplectic form) to $\pi^*(dz_1)$. It preserves the fibres of the Lagrangian fibration over $\Delta$ and a characteristic $1$-cycle is essentially an integral curve of $v_1$. Similarly, there are vector fields $v_2,\ldots,v_n$ dual to $\pi^*(dz_2),\ldots,\pi^*(dz_n)$ which also preserve the fibres, and generate a $\C^{n-1}$-action. The singular locus of a fibre will be an orbit of this $\C^{n-1}$-action. For the compactified Jacobian $\overline{J}$ of a cuspidal curve, the characteristic $1$-cycle would be {\em tangent\/} to the singular locus, implying that $v_1$ lies in the span of $v_2,\ldots,v_n$. This contradicts the independence of $dz_1,dz_2,\ldots,dz_n$.
\end{rmk}

`Cuspidal' fibres can appear in codimension one in Lagrangian fibrations, but only if the characteristic $1$-cycles are themselves cuspidal curves.

\begin{exm}
Let $S\rightarrow\P^1$ be an elliptic K3 surface containing a cuspidal rational curve $E_p$ over the point $p\in\P^1$. If $p_2,\ldots,p_n$ are generic points of $\P^1$, distinct from $p$ and each other, then the fibre of $\mathrm{Hilb}^nS\rightarrow\mathrm{Sym}^n\P^1$ over $\{p,p_2,\ldots,p_n\}$ is $E_p\times E_2\times\ldots\times E_n$. This is the kind of `cuspidal' fibre allowed by Hwang and Oguiso's classification.
\end{exm}

\begin{lem}
\label{d=1}
The compactified relative Jacobian $X=\overline{J}^d(Y/B)$ is a Lagrangian fibration if and only if $X^1=\overline{J}^1(Y/B)$ is a Lagrangian fibration.
\end{lem}

\begin{prf}
The following argument is inspired by Section~2 of Markushevich~\cite{markushevich96}, particularly Proposition~2.3. Without loss of generality, suppose that $X\rightarrow B$ is a Lagrangian fibration, with holomorphic symplectic form $\sigma$ on $X$. Because the curves in the family $Y/B$ are integral, we can choose local sections $s_i:U_i\rightarrow Y|_{U_i}$, where $\{U_i\}$ is an open cover of $B$ (in the analytic topology). Then tensoring with $\O_Y((d-1)s_i)$ gives a local isomorphism
$$\phi_i:X^1|_{U_i}\rightarrow X|_{U_i}$$
which commutes with the projections to $U_i$. In other words, $X$ and $X^1$ are locally isomorphic as fibrations. Define a two-form
$$\sigma_i:=\phi_i^*\left(\sigma|_{X|_{U_i}}\right)$$
on $X^1|_{U_i}$. We want to glue these two-forms together to get a holomorphic symplectic form on $X^1$; however, they may not agree on the overlaps $X^1|_{U_i\cap U_j}$.

Restrict to the locus of non-critical points of the fibration $\pi:X^1\rightarrow B$. If $p$ is a non-critical point, then we can pull-back coordinates from $B$ and complete them to coordinates $(x_1,\ldots,x_n,y_1,\ldots,y_n)$ around $p$, with $(y_1,\ldots,y_n)$ in the fibre direction. These coordinates can be chosen so that
$$\sigma_i=\sum_k dx_k\wedge dy_k + \sum_{k,l}f^i_{kl}dx_k\wedge dx_l$$
and
$$\sigma_j=\sum_k dx_k\wedge dy_k + \sum_{k,l}f^j_{kl}dx_k\wedge dx_l,$$
where $f^i_{kl}$ and $f^j_{kl}$ are local functions on $X^1$. Note that there are no $dy_k\wedge dy_l$ terms because $T_pX^1_t\subset T_pX^1$ is mapped to the Lagrangian subspace $T_{\phi_i(p)}X_t\subset T_{\phi_i(p)}X$ by $d\phi_i$ (and similarly for $\phi_j$). The difference is
$$\sigma_i-\sigma_j=\sum_{k,l}g^{ij}_{kl}dx_k\wedge dx_l.$$
Now the fibre coordinates $(y_1,\ldots,y_n)$ extend to global coordinates, modulo a lattice, on the smooth fibres. It follows that the holomorphic functions $g^{ij}_{kl}$ are independent of $(y_1,\ldots,y_n)$ on the smooth fibres, and therefore everywhere by continuity. Since the coordinates $(x_1,\ldots,x_n)$ are pulled back from the base $B$, we see that
$$\sigma_i-\sigma_j=\pi^*\beta_{ij}$$
where $\beta_{ij}$ is a two-form on $U_i\cap U_j$. Moreover, $[\beta_{ij}]$ defines a class in $\H^1(B,\Omega^2_B)$. But $B\cong\P^n$ by Lemma~\ref{base}, and so this cohomology group vanishes. Therefore we can write $\beta_{ij}=b_i-b_j$, where $b_i$ and $b_j$ are two-forms on $U_i$ and $U_j$, respectively. The modified two-forms
$$\sigma_i-\pi^*b_i\qquad\mbox{and}\qquad\sigma_j-\pi^*b_j$$
agree on $X^1|_{U_i\cap U_j}$ and can be glued to give a two-form $\sigma^1$ on $X^1$, or at least on its non-critical locus.

It remains to show that $\sigma^1$ extends to all of $X^1$ and is non-degenerate. The critical locus of the Lagrangian fibration $\pi:X\rightarrow B$ has codimension at least two. Since $X$ and $X^1$ are locally isomorphic as fibrations, the critical locus of $\pi:X^1\rightarrow B$ will also have codimension at least two. By Hartog's Theorem, $\sigma^1$ extends to a holomorphic two-form on all of $X^1$. Finally, the degeneracy locus of $\sigma^1$ is the same as the zero locus of $\wedge^n\sigma^1\in\H^0(X^1,K_{X^1})$, i.e., it has codimension one if it is non-empty. But the local form
$$\sigma_i-\pi^*b_i=\sum_k dx_k\wedge dy_k + \sum_{k,l}h^i_{kl}dx_k\wedge dx_l$$
of $\sigma^1$ shows that it is non-degenerate on the non-critical locus of $X^1$, whose complement has codimension at least two. Therefore $\sigma^1$ is non-degenerate everywhere, and $X^1\rightarrow B$ is a Lagrangian fibration.
\end{prf}

\begin{lem}
\label{locallyfree}
Use $\pi$ to denote the projections of both $Y$ and $X$ to $B\cong\P^n$.
\begin{enumerate}
\item The first direct image sheaf $R^1\pi_*\mathcal{O}_X$ is isomorphic to $\Omega^1_{\P^n}$. 
\item The first direct image sheaf $R^1\pi_*\mathcal{O}_Y$ is isomorphic to $\Omega^1_{\P^n}$. 
\end{enumerate}
\end{lem}

\begin{prf}
The first statement of the lemma was proved by Matsushita~\cite{matsushita05}; it is true for any Lagrangian fibration $\pi:X\rightarrow\P^n$.

For the second statement, first notice that by Lemma~\ref{genusn} the arithmetic genus $\mathrm{dim}\mathrm{H}^1(C,\mathcal{O}_C)=n$ for all curves $C$ in the family $Y/B$. Therefore the dimension of the fibres of
$R^1\pi_*\mathcal{O}_Y$ does not jump, and $R^1\pi_*\mathcal{O}_Y$ is locally free of rank $n$. 

By Lemma~\ref{d=1}, we may assume that the degree $d=1$. Given an integral curve $C$, its Abel map $C\hookrightarrow\overline{J}$ is a closed embedding, and canonical for $d=1$. The relative Abel map $Y\hookrightarrow X$ is also a canonical closed embedding (see Altman and Kleiman, Theorem~8.8~\cite{ak80}). We identify $Y$ with its image in $X$. The short exact sequence
$$0\rightarrow\mathcal{I}_Y\rightarrow\mathcal{O}_X\rightarrow\mathcal{O}_Y\rightarrow
0,$$
where $\mathcal{I}_Y$ is the ideal sheaf of $Y\subset X$, yields the
long exact sequence
$$R^1\pi_*\mathcal{I}_Y\rightarrow
R^1\pi_*\mathcal{O}_X\stackrel{\alpha}{\longrightarrow}R^1\pi_*\mathcal{O}_Y\rightarrow
R^2\pi_*\mathcal{I}_Y.$$
We must show that $\alpha$ is an isomorphism. Since $R^1\pi_*\mathcal{O}_X$
and $R^1\pi_*\mathcal{O}_Y$ are both locally free of rank $n$, it suffices to show
that $\alpha$ is an isomorphism on the complement of a codimension two
subset of $\P^n$.

Above a generic point $t\in\P^n$, the fibre of $Y$ is a smooth curve
$Y_t$ and the fibre $X_t$ is the Jacobian of $Y_t$. In this case $\alpha_t$ is a composition of canonical isomorphisms
$$(R^1\pi_*\mathcal{O}_Y)_t\cong\mathrm{H}^1(Y_t,\mathcal{O}_{Y_t})\cong\mathrm{H}^1(X_t,\mathcal{O}_{X_t})\cong(R^1\pi_*\mathcal{O}_X)_t.$$
Now let $t$ be a generic point of the discriminant locus
$\Delta\subset\P^n$. By Lemma~\ref{singlenode} the curve $C=Y_t$ has a single node, and its compactified Jacobian $\overline{J}=X_t$ is as described in the previous subsection: recall that the normalization $P$ of $\overline{J}$ is a $\P^1$-bundle over the Jacobian $\tilde{J}$ of the normalization $\tilde{C}$ of $C$. Let $g$ denote both the normalization maps $\tilde{C}\rightarrow C$ and $P\rightarrow\overline{J}$. We obtain the following commutative
diagram involving short exact sequences
$$\begin{array}{ccccccccc}
0 & \rightarrow & \mathcal{O}_{\overline{J}} & \rightarrow &
g_*\mathcal{O}_P & \rightarrow & \mathcal{G} & \rightarrow
& 0 \\
  & & \downarrow & & \downarrow & & \downarrow & & \\
0 & \rightarrow & \mathcal{O}_C & \rightarrow &
g_*\mathcal{O}_{\tilde{C}} & \rightarrow & \mathcal{G}^{\prime} & \rightarrow
& 0 \\
\end{array}$$
where $\mathcal{G}$ is supported on the singular locus of $\overline{J}$ and
$\mathcal{G}^{\prime}$ is supported on the node of $C$. The vertical
arrows come from the Abel embedding of $C$ in $\overline{J}$ and the
corresponding embedding of $\tilde{C}$ in $P$. Taking
cohomology we get
$$\begin{array}{cccccccc}
\stackrel{0}{\longrightarrow} & {\H}^0(\overline{J},\mathcal{G}) & \rightarrow &
{\H}^1(\overline{J},\mathcal{O}_{\overline{J}}) & \stackrel{i_1}{\longrightarrow} &
{\H}^1(\overline{J},g_*\mathcal{O}_P)\cong{\H}^1(P,\mathcal{O}_P)
& \rightarrow & \ldots \\
 & \downarrow & & \downarrow & & \downarrow & & \\
\stackrel{0}{\longrightarrow} & {\H}^0(C,\mathcal{G}^{\prime}) & \rightarrow &
{\H}^1(C,\mathcal{O}_C) & \stackrel{i_2}{\longrightarrow} &
{\H}^1(C,g_*\mathcal{O}_{\tilde{C}})\cong{\H}^1(\tilde{C},\mathcal{O}_{\tilde{C}})
& \rightarrow & \ldots \\
\end{array}$$
The first term in each row is isomorphic to $\C$, and the left-most
vertical arrow is an isomorphism. Since $P$ is a
$\P^1$-bundle over $\tilde{J}$, we have
$$\H^1(P,\mathcal{O}_P)\cong\H^1(\tilde{J},\mathcal{O}_{\tilde{J}})\cong\H^1(\tilde{C},\mathcal{O}_{\tilde{C}}).$$
Therefore the right-most vertical arrow is also an
isomorphism. Since $\H^1(\tilde{C},\mathcal{O}_{\tilde{C}})$ has dimension $n-1$ and $\H^1(C,\mathcal{O}_C)$ has dimension $n$, the map labelled $i_2$ must be surjective. Similarly, $\H^1(P,\mathcal{O}_P)$ has dimension $n-1$ and $\H^1(\overline{J},\mathcal{O}_{\overline{J}})$ has dimension at least $n$ (by semi-continuity applied to the family $X\rightarrow\P^n$), so the map labelled $i_1$ must also be surjective. Therefore the middle vertical arrow must be an isomorphism, and again $\alpha_t$ is a composition of isomorphisms
$$(R^1\pi_*\mathcal{O}_Y)_t\cong\H^1(C,\mathcal{O}_C)\cong\H^1(\overline{J},\mathcal{O}_{\overline{J}})\cong(R^1\pi_*\mathcal{O}_X)_t$$
for a generic point $t\in\Delta$. This completes the proof.
\end{prf}

\subsection{The Beauville-Mukai system}

Let $S$ be a K3 surface containing a smooth genus $n$ curve $C$. Then $C$ moves in an $n$-dimensional linear system $|C|\cong\P^n$; let $Y\rightarrow\P^n$ be the family of curves linearly equivalent to $C$. If $S$ is generic in the sense that $\mathrm{Pic}S\cong\Z C$, then every curve in this family will be integral. 

\begin{dfn}
The compactified relative Jacobian $X=\overline{J}^d(Y/\P^n)$ is known as the Beauville-Mukai integrable system~\cite{beauville99, mukai84}.
\end{dfn}
By thinking of a rank-one torsion free sheaf $\mathcal{F}$ on a curve $C\stackrel{\iota}{\hookrightarrow}S$ as a torsion sheaf $\iota_*\mathcal{F}$ on $S$, one can identify $X$ with an irreducible component of the Mukai moduli space of stable sheaves on $S$. From this it follows that
\begin{itemize}
\item $X$ is smooth of dimension $2n$,
\item $X$ admits a holomorphic symplectic structure,
\item $X$ is deformation equivalent to $\mathrm{Hilb}^nS$.
\end{itemize}
General results of Matsushita~\cite{matsushita99, matsushita00} then imply that the Beauville-Mukai system $X\rightarrow\P^n$ is a Lagrangian fibration.

\begin{rmk}
One can still construct a Beauville-Mukai system when the linear system $|C|$ contains reducible and/or non-reduced curves, by {\em defining\/} the compactified relative Jacobian $X$ to be the corresponding Mukai moduli space, as above. However, $X$ may be non-compact in some cases (depending on $d$); adding semi-stable sheaves compactifies the moduli space, but also introduces singularities.
\end{rmk}

\section{Statement of the theorem}

By Lemma~\ref{base} we can assume that the base $B$ of our family of curves is $\P^n$. The main goal of this paper is to prove the following theorem, whose proof will take up all of Section~4.

\begin{thm}
\label{main}
Let $Y\rightarrow\P^n$ be a flat family of integral Gorenstein curves whose compactified relative Jacobian $X=\overline{J}^d(Y/\P^n)$ is a Lagrangian fibration. Then
\begin{enumerate}
\item the degree of the discriminant locus $\Delta\subset\P^n$ is at least $4n+2$,
\item if $\mathrm{deg}\Delta>4n+20$ then $X$ is a Beauville-Mukai integrable system, i.e., the family of curves is a complete linear system of curves in a K3 surface. 
\end{enumerate}
\end{thm}

\begin{rmk}
The discriminant locus $\Delta\subset\P^n$ of the Beauville-Mukai system has degree $6n+18$, as calculated in Section~5 of~\cite{sawon08i}.
\end{rmk}


\begin{rmk}
The hypothesis that the curves be Gorenstein is needed only to prove Lemma~\ref{Ysmooth} below; moreover, we expect that it follows from the other hypotheses.

Let $t\in\P^n$ be an arbitrary point and let $H_1,\ldots,H_n$ be $n$ hyperplanes in $\P^n$ meeting transversally at $t$. Then the compactified Jacobian $\overline{J}^d=X_t$ of $C=Y_t$ is a complete intersection 
$$\pi^{-1}(H_1)\cap\ldots\pi^{-1}(H_n),$$
and in particular it is Gorenstein. We expect that the curve $C$ must automatically be Gorenstein in this case: If $C$ has a non-Gorenstein singularity at $p$ then its Hilbert scheme $\mathrm{Hilb}^m(C)$ will have a non-Gorenstein singularity at $\{p,p_2,\ldots,p_m\}$, where $p,p_2,\ldots,p_m$ are distinct point with $p_2,\ldots,p_m$ lying in the smooth locus of $C$. Now we invoke the Abel map
$$\mathrm{Hilb}^m(C)\rightarrow\overline{J}^m\cong\overline{J}^d$$
to compare singularities of $\mathrm{Hilb}^m(C)$ to those of $\overline{J}^d$. Unfortunately, D'Souza and Rego (see Altman and Kleiman, Theorem~8.6~\cite{ak80}) showed that the Abel map is smooth if and only if $m\geq 2n-1$ {\em and }$C$ {\em is Gorenstein\/}, so this argument is unfortunately circular.
\end{rmk}

As a corollary of Theorem~\ref{main}, we obtain the following result,
first proved by Markushevich in~\cite{markushevich96}.

\begin{cor}
\label{corollary}
Let $Y\rightarrow\P^2$ be a flat family of integral Gorenstein curves whose compactified relative Jacobian $X=\overline{J}^d(Y/\P^2)$ is a Lagrangian fibration. Then $X$ is a Beauville-Mukai integrable system.
\end{cor}

\begin{prf}
It suffices to prove that the hypothesis $\mathrm{deg}\Delta >
4n+20$ in Theorem~\ref{main} is automatically satisfied when $n=2$.
In~\cite{sawon08i} it was proved that if $X\rightarrow\P^n$
is a Lagrangian fibration by principally polarized abelian varieties
with `good singular fibres' then
$$\mathrm{deg}\Delta=24\left(n!\sqrt{\hat{A}}[X]\right)^{\frac{1}{n}}$$
where $\sqrt{\hat{A}}[X]$ is the characteristic number of $X$ given by
the square root of the $\hat{A}$-polynomial. For our family $Y\rightarrow\P^2$, the singular curves in codimension one have single nodes by Lemma~\ref{singlenode}, which ensures that $X\rightarrow\P^2$ has good singular fibres. It remains to bound $\sqrt{\hat{A}}[X]$ from below.

We claim that an irreducible holomorphic symplectic fourfold satisfies
$$\sqrt{\hat{A}}[X]\geq \frac{25}{32}.$$
As in Section 8 of~\cite{sawon08i}, we can write $\sqrt{\hat{A}}[X]$ in
terms of Betti numbers and then use Guan's bounds~\cite{guan01} on
these Betti numbers. The relation
$$\hat{A}[X]=\frac{1}{720}(3c_2^2[X]-c_4[X])=\chi(\O_X)=3$$
between the Chern numbers allows us to write $\sqrt{\hat{A}}[X]$
solely in terms of $c_4[X]$, giving
$$\sqrt{\hat{A}}[X]=\frac{3024-c_4[X]}{3456}.$$
Since $X$ is simply connected, the first Betti number vanishes. The
fourth Betti number is determined by Salamon's relation~\cite{salamon96}
$$b_4=46+10b_2-b_3,$$
and therefore the topological Euler characteristic of $X$ is
$$c_4[X]=48+12b_2-3b_3.$$
This gives
$$\sqrt{\hat{A}}[X]=\frac{992-4b_2+b_3}{1152}.$$
Now Guan~\cite{guan01} proved that $b_2$ and $b_3$ can take only
finitely many values; the smallest value for $\sqrt{\hat{A}}[X]$,
$\frac{25}{32}$, occurs when $b_2=23$ and $b_3=0$. This establishes
the claim.

Finally, when $n=2$
$$\mathrm{deg}\Delta=24\left(2!\sqrt{\hat{A}}[X]\right)^{\frac{1}{2}}\geq
24\left(2!\frac{25}{32}\right)^{\frac{1}{2}}=30,$$
which is greater than $4n+20$.
\end{prf}

\section{The proof of Theorem 8}

\subsection{Outline}

Let us begin with an outline of the proof. We first show that the
total space of the family of curves $Y\rightarrow\P^n$ is
smooth. By Lemma~\ref{d=1}, we may assume that the degree $d=1$, so that the relative Abel map gives a canonical closed embedding $Y\hookrightarrow X$. Henceforth we identify $Y$ with its image in $X$.

Let $\sigma$ be the holomorphic symplectic form on $X$. The restriction $\sigma|_Y\in\H^0(Y,\Omega^2_Y)$ must be degenerate, and the kernel of $\sigma|_Y$ defines a distribution $F\subset TY$. More precisely, $F$ is the kernel of the morphism of sheaves
\begin{eqnarray*}
TY & \rightarrow & \Omega^1_Y \\
v & \mapsto & \sigma|_Y(v,-)
\end{eqnarray*}
where $v$ is a local section of $TY$. As a coherent subsheaf of the locally free sheaf $TY$, $F$ will be locally free on the complement of a codimension two subset. A priori $F$ need not be locally free on all of $Y$, because the rank of $\sigma|_Y$ need not
be constant on $Y$. We do know that the rank of $\sigma|_{T_yY}$ is even for each $y\in Y$ and that
$$2n-2\mathrm{codim}Y=2\leq\mathrm{rank}\sigma|_{T_yY}\leq
2\left\lfloor\frac{\mathrm{dim}Y}{2}\right\rfloor=2\left\lfloor\frac{n+1}{2}\right\rfloor.$$
We also know that the rank is semi-continuous, meaning that
$$\{y\in Y|\mathrm{rank}\sigma|_{T_yY}\leq m\}$$
is a closed subset of $Y$ for each $m$. In the course of the proof we
will prove that $\sigma|_Y$ has constant rank $2$, and hence $F$ will
be locally free of rank $n-1$.

Since $\sigma$ is $d$-closed the distribution is integrable, i.e., if
$u$ and $v\in F$ then
$$\sigma|_Y([u,v],w)=d\sigma|_Y(u,v,w)+u\sigma|_Y(v,w)-v\sigma|_Y(u,w)=0$$
so $[u,v]\in F$. It therefore defines a foliation, known as the {\em characteristic foliation\/}. This notion of
integrability (of the distribution) should not be confused with the
following notion of algebraic integrability (of the foliation).

We will show that the characteristic foliation is algebraically integrable, i.e.,
that the leaves are algebraic and in particular compact. This is the
most difficult part of the proof, and our approach is inspired by the
following theorem of Kebekus, Sol{\'a} Conde, and Toma~\cite{ksct07}.
\begin{thm}
\label{kebekus}
Let $Y$ be a smooth variety, $C$ a complete curve in $Y$ and $F\subset
TY$ a foliation which is regular over $C$. Assume that $F|_C$ is an
ample vector bundle. Then the leaf through any point of $C$ is
algebraic, and hence its closure is compact. Moreover, the closure of
the leaf through the generic point of $C$ is rationally connected. If
in addition $F$ is regular everywhere then all leaves are rationally
connected.
\end{thm}
Once we have shown algebraic integrability, the space of leaves $Y/F$
will be a well-defined surface $S$. Moreover, the two-form $\sigma|_Y$
will descend to a non-degenerate two-form on $S$ under the projection
$Y\rightarrow S$, implying that $S$ is either a K3 surface or an
abelian surface. We will show that each curve in the family
$Y\rightarrow\P^n$ maps isomorphically to its image under the
projection $Y\rightarrow S$. Therefore $S$ contains an $n$-dimensional
linear system of genus $n$ curves, and it follows that $S$ must be a
K3 surface. This will conclude the proof.

\subsection{Rational curves in $Y$}

Our goal in this subsection is to show that $Y$ contains (many) rational curves.

\begin{lem}
\label{Ysmooth}
The total space $Y$ of the family of curves is smooth.
\end{lem}

\begin{prf}
Suppose that $Y$ is singular at a point $p$. Let $t=\pi(p)\in\P^n$, and let $p_2,\ldots,p_m$ be distinct points in the smooth locus of the curve $Y_t$, which are also distinct from $p$. Then the Hilbert scheme $\mathrm{Hilb}^m(Y/\P^n)$ will have a singularity at $\{p,p_2,\ldots,p_m\}$. By hypothesis, the curves in the family $Y/\P^n$ are Gorenstein, so for $m\geq 2n-1$ the Abel map
$$\mathrm{Hilb}^m(Y/\P^n)\rightarrow\overline{J}^m(Y/\P^n)$$
is smooth by a theorem of D'Souza and Rego (see Altman and Kleiman, Theorem~8.6~\cite{ak80}). Therefore $\overline{J}^m(Y/\P^n)$ will be singular at the image of $\{p,p_2,\ldots,p_m\}$. Now $\overline{J}^m(Y/\P^n)$ is locally isomorphic (as a fibration over $\P^n$) to $X=\overline{J}^1(Y/\P^n)$, and this means that $X$ must also be singular. But this contradicts the hypothesis that $X$ is a holomorphic symplectic manifold, and in particular smooth.
\end{prf}

\begin{rmk}
The above lemma is the only place where we require the curves in the family $Y\rightarrow\P^n$ to be Gorenstein; indeed, the Gorenstein hypothesis could be replaced by the hypothesis that the total space $Y$ be smooth.
\end{rmk}

Let $t\in\P^n$ and let $p$ be a point in the smooth locus of the curve $Y_t$. Let $X_t$ be the corresponding fibre of $\pi:X\rightarrow\P^n$. There is a short exact sequence
$$0\rightarrow T_pX_t\rightarrow T_pX \rightarrow\pi^*T_t\P^n\rightarrow 0$$
and $T_pY_t$ sits inside the first term since $Y_t$ sits inside its compactified Jacobian $X_t$. Now a non-zero vector
$$w\in T_pY_t\subset T_pX_t\subset T_pX$$
defines a linear functional $\sigma_p(w,-)$ on $T_pX$, where $\sigma$ is the holomorphic symplectic form on $X$. This in turn defines a linear functional on $T_t\P^n$, which takes $v$ in $T_t\P^n\cong \pi^*T_t\P^n$ to $\sigma_p(w,v^{\dagger})$, where $v^{\dagger}$ is a lift of $v$ to $T_pX$. Note that $\sigma_p(w,v^{\dagger})$ is independent of the choice of lift $v^{\dagger}$ of $v$ since two lifts will differ by an element of $T_pX_t$, and $T_pX_t\subset T_pX$ is a Lagrangian subspace. This linear functional on $T_t\P^n$ will vanish on a hyperplane $H_p$ in $T_t\P^n$.

\begin{dfn}
We call $H_p$ the associated hyperplane of $p$. It can be regarded as a hyperplane in $T_t\P^n$ or as a hyperplane in $\P^n$ passing through $t$.
\end{dfn}

Let $t$ be a generic point of $\P^n$, and let $\ell\subset\P^n$ be a line through $t$. If $\ell$ is chosen generically, then $Z=\pi^{-1}(\ell)$ is a smooth surface which is fibred by genus $n$ curves over $\ell\cong\P^1$ (here $\pi$ denotes the projection $Y\rightarrow\P^n$, not $X\rightarrow\P^n$).

\begin{lem}
\label{nonvanishing}
Let $p$ be a point in the smooth locus of the curve $Y_t$. If $\ell$ does not lie in the associated hyperplane $H_p$ of $p$, then the restriction $\sigma|_Z$ of the holomorphic symplectic form $\sigma$ to $Z$ does
not vanish at $p$. In particular, for generic $\ell$, $\sigma|_Z$ is a non-trivial holomorphic two-form on $Z$ and hence $h^{2,0}(Z)=h^{0,2}(Z)> 0$.
\end{lem}

\begin{prf}
Suppose that $\ell$ is not contained in $H_p$. Let $w\in T_pY_t\subset T_pX$ and
$$v\in T_t\ell\subset T_t\P^n\cong\pi^*T_t\P^n$$
be non-zero vectors. Then $w$ lies in $T_pZ$ and we can lift $v$ to $v^{\dagger}$ in $T_pZ\subset T_pX$. Since $v\not\in H_p$, $\sigma_p(w,v^{\dagger})\neq 0$, and thus $\sigma|_Z$ does not vanish at $p$.
\end{prf}

\begin{lem}
\label{cohomology}
For generic $\ell$, the cohomology of $\mathcal{O}_Z$ is given by
$$h^{0,k}(Z)=\left\{\begin{array}{cl}
                 1 & k=0, \\
                 0 & k=1, \\
                 1 & k=2. \\
                 \end{array}\right.$$
\end{lem}

\begin{prf}
By Lemma~\ref{locallyfree} we have
$R^1\pi_*\mathcal{O}_Y\cong\Omega^1_{\P^n}$, and therefore
$$R^1\pi_*\mathcal{O}_Z\cong R^1\pi_*\mathcal{O}_Y|_{\ell}\cong
\mathcal{O}_{\P^1}(-2)\oplus\mathcal{O}_{\P^1}(-1)^{\oplus(n-1)}.$$
Inserting this into the Leray spectral sequence
$$\mathrm{H}^i(\P^1,R^j\pi_*\mathcal{O}_Z)\Rightarrow\mathrm{H}^{i+j}(Z,\mathcal{O}_Z)$$
for $\pi:Z\rightarrow\ell$ yields the result.
\end{prf}

We can now prove the first statement of Theorem~\ref{main}.

\begin{lem}
The degree of the discriminant locus $\Delta\subset\P^n$ is at least $4n+2$.
\end{lem}

\begin{prf}
By the previous lemma $Z$ has $q=0$ and $p_g=1$, so it must have a unique minimal model of Kodaira dimension $0$, $1$, or $2$. If $h^{1,1}(Z)=1$ then $Z$ itself must be minimal. Moreover, $Z$ will have Euler characteristic $e(Z)=5$ and $K_Z^2=19$ (by Noether's formula). Therefore $Z$ is of general type; but then $K_Z^2>3e(Z)$ contradicts the Bogomolov-Miyaoka-Yau inequality. We conclude that $h^{1,1}\geq 2$, and thus $e(Z)\geq 6$.

Since $\ell$ is generic, it will intersect $\Delta$ in $\mathrm{deg}\Delta$
distinct points, and if $t\in\ell\cap\Delta$ then the curve $Y_t$ will
contain precisely one node. A fibration by genus $n$ curves over
$\P^1$ with $\mathrm{deg}\Delta$ nodal fibres will have Euler
characteristic
$$e(Z)=-2(2n-2)+\mathrm{deg}\Delta.$$
Therefore
$$\mathrm{deg}\Delta = e(Z)+4n-4\geq 4n+2.$$
\end{prf}

We now turn to the second statement of Theorem~\ref{main}. For the remainder of this section we assume that $\mathrm{deg}\Delta>4n+20$.

\begin{lem}
\label{notminimal}
For generic $\ell$, the surface $Z$ is not a minimal surface: it contains at least one
$(-1)$-curve, i.e., rational curve with self-intersection $-1$.
\end{lem}

\begin{prf}
In the previous proof we calculated
$$e(Z)=-2(2n-2)+\mathrm{deg}\Delta.$$
Combining this with Noether's formula
$$K_Z^2+e(Z)=12\chi(\mathcal{O}_Z)=24$$
gives
$$K_Z^2=24-e(Z)=4n+20-\mathrm{deg}\Delta,$$
which is negative by the hypothesis $\mathrm{deg}\Delta >4n+20$.

Suppose that $Z$ is minimal. Since $K_Z^2<0$ it must have Kodaira
dimension $-\infty$; but then $p_g$ would have to be zero,
contradicting the existence of a non-trivial holomorphic two-form
$\sigma|_Z$ on $Z$.
\end{prf}

\begin{rmk}
The hypothesis $\mathrm{deg}\Delta >4n+20$ is necessary to show
that $K_Z^2<0$. Without this, the Kodaira dimension of $Z$ could be
zero, one, or two. In the Appendix we discuss minimal
surfaces $Z$ with $\mathrm{kod}(Z)\in\{0,1,2\}$, $q=0$, and
$p_q=1$. Moreover, we consider whether such surfaces can be fibred
over $\P^1$ by genus $n$ curves.
\end{rmk}

The following two lemmas show that the $(-1)$-curves we have found sit fairly generically inside $Y$, which will be important later.

\begin{lem}
\label{notcontained}
Let $p$ be a point in the smooth locus of the curve $Y_t$. If $\ell$ does not lie in the associated hyperplane $H_p$ of $p$, then no $(-1)$-curve
in $Z$ contains $p$.
\end{lem}

\begin{prf}
Since $\sigma|_Z$ is a section of the canonical bundle
$K_Z=\Omega^2_Z$, it must vanish on every $(-1)$-curve. However, if $\ell\not\subset H_p$ then by Lemma~\ref{nonvanishing} $\sigma|_Z$ does not vanish at $p$.
\end{prf}

With a little more effort, we can extend the above lemma to nodes $p$ of $Y_t$.

\begin{lem}
\label{notcontained2}
Let $t$ be a generic point of $\Delta\subset\P^n$, so that $Y_t$ contains a single node $p$. If $\ell$ is a generic line through $t$ then no $(-1)$-curve in $Z$ contains $p$.
\end{lem}

\begin{prf}
A local calculation shows that the surface $Z=\pi^{-1}(\ell)$ will have a singularity at $p$ if $\ell$ is tangent to $\Delta$ at $t$. However, we are assuming $\ell$ is a generic line through $t$, so the surface $Z$ will be smooth. The structure of the compactified Jacobian $X_t$ of $Y_t$ is described in the first example of Subsection~2.1. Let
$$w\mbox{ and }w^{\prime}\in T_pY_t\subset T_pX_t\subset T_pX$$
be non-zero vectors tangent to the two branches of $Y_t$ at $p$. We can complete these to bases $\{w,w_2,\ldots,w_n\}$ and $\{w^{\prime},w_2,\ldots,w_n\}$ for the tangent spaces to the two branches of $X_t$ at $p$. Then
$$\mathrm{span}\{w,w_2,\ldots,w_n\}\qquad\mbox{and}\qquad\mathrm{span}\{w^{\prime},w_2,\ldots,w_n\}$$
will both be Lagrangian subspaces of $T_pX$. Suppose that $p$ lies in a $(-1)$-curve in $Z$. Then $\sigma|_Z$ vanishes at $p$, and since $w$ and $w^{\prime}$ both lie in $T_pY_t\subset T_pZ$, we have
$$\sigma(w,w^{\prime})=0.$$
But then the $(n+1)$-dimensional subspace
$$\mathrm{span}\{w,w^{\prime},w_2,\ldots,w_n\}\subset T_pX$$
is isotropic, contradicting the non-degeneracy of $\sigma$ on $T_pX$.
\end{prf}

\subsection{Rationality of the leaves}

We wish to use Theorem~\ref{kebekus}, due to Kebekus et al., to study
the characteristic foliation $F$ on $Y$. An obvious curve to use would be a fibre
$Y_t$ of $Y\rightarrow\P^n$; however, one can easily show that
$F|_{Y_t}$ will not be ample. Instead, we will eventually show that the rational $(-1)$-curves found in the previous subsection lie inside the leaves of $F$. But first we show that the leaves of the foliation meet  the curve $Y_t$ transversally at smooth points of $Y_t$. 

\begin{lem}
\label{transverse}
Let $t\in\P^n$ and let $p$ be a point in the smooth locus of the curve $Y_t$. A non-zero vector $w\in
T_pY_t\subset T_pY$ cannot lie in $F\subset TY$.
\end{lem}

\begin{prf}
Suppose conversely that $w\in F_p$. This means that $\sigma_p(w,v)=0$
for every vector $v\in T_pY$. We have the following
commutative diagram:
$$\begin{array}{ccccccccc}
  0 & \rightarrow & T_pX_t & \rightarrow & T_pX & \rightarrow &
  \pi^*T_t\P^n & \rightarrow & 0 \\
    & & \cup & & \cup & & \| & & \\
  0 & \rightarrow & T_pY_t & \rightarrow & T_pY & \rightarrow &
  \pi^*T_t\P^n & \rightarrow & 0 \\
\end{array}$$
The map from $T_pX$ to $\pi^*T_t\P^n$ is the differential $d\pi_p$. Let $u$ be an arbitrary vector in $T_pX$ and let $d\pi_p(u)$ be its image in $\pi^*T_t\P^n$. We can lift $d\pi_p(u)$ to $v\in T_pY\subset T_pX$. Then $u$ and $v$
will both be lifts of $d\pi_p(u)$, so they will differ only by a
vector in $T_pX_t$. But $T_pX_t$ is a
Lagrangian subspace of $T_pX$, and $w\in T_pY_t\subset T_pX_t$, so
$$\sigma_p(w,u-v)=0.$$
We therefore have
$$\sigma_p(w,u)=\sigma_p(w,v)=0.$$
Since this is true for all $u\in T_pX$, the holomorphic symplectic form $\sigma$ on $X$
will be degenerate at $p$, a contradiction.
\end{prf}

\begin{lem}
\label{constantrank}
The distribution $F\subset TY$ is locally free of rank $n-1$. Equivalently, the
rank of $\sigma|_Y$ is constant and equal to two.
\end{lem}

\begin{prf}
Suppose that $\sigma|_Y$ has rank $2m$ at a generic point. By
semi-continuity $\sigma|_Y$ drops rank on closed subsets
of $Y$. We will show that $2m=2$, and therefore
$\mathrm{rank}\sigma|_Y$ is constant because the rank cannot drop
below two.

The kernel $F$ of $\sigma|_Y$ has rank $n+1-2m$ at a generic
point. Since $F$ is a coherent subsheaf of the locally free sheaf
$TY$, it will be locally free on the complement of a codimension two
subset. Since $F$ is the kernel of $TY\rightarrow\Omega^1_Y$, we must
have an injection $TY/F\rightarrow\Omega^1_Y$ by the universal
property of the quotient. Thus $TY/F$ is also locally free on the
complement of a codimension two subset, since it is a coherent
subsheaf of the locally free sheaf $\Omega^1_Y$. So suppose that both
$F$ and $TY/F$ are locally free on the complement of the codimension
two subset $\Sigma\subset Y$. Now $TY/F$ has rank $2m$ at a generic
point, and $\sigma|_Y$ descends to a two-form $\omega$ on $TY/F$ which
is non-degenerate at a generic point; in particular, $\omega^{\wedge m}$
is a non-trivial section of $\Lambda^{2m}(TY/F)^{\vee}$. We illustrate
this with two examples before continuing.

\begin{exm}
Suppose that $n=4$ (so that $Y$ is $5$-dimensional), $2m=4$, and
$\sigma|_Y$ looks locally like
$$dz_1\wedge dz_2+fdz_3\wedge dz_4$$
where $f$ is a function on $Y$. Although $\sigma|_Y$ drops rank along
the hypersurface $\{f=0\}$, the vectors
$\frac{\partial\phantom{z}}{\partial z_3}$ and
$\frac{\partial\phantom{z}}{\partial z_4}$ do not extend to local
vector fields in the kernel of the map $TY\rightarrow\Omega^1_Y$. So
$F$ is simply generated by $\frac{\partial\phantom{z}}{\partial
  z_5}$. The quotient $TY/F$ is generated by
$$\left\{\frac{\partial\phantom{z}}{\partial
  z_1},\frac{\partial\phantom{z}}{\partial
  z_2},\frac{\partial\phantom{z}}{\partial
  z_3},\frac{\partial\phantom{z}}{\partial z_4}\right\}$$
(where we have identified $\frac{\partial\phantom{z}}{\partial z_i}$
  with its coset in $TY/F$), $\omega$ looks just like $\sigma|_Y$, and
$$\omega^{\wedge 2}=fdz_1\wedge dz_2\wedge dz_3\wedge dz_4$$
is a non-trivial section of $\Lambda^4(TY/F)^{\vee}$, which vanishes on $\{f=0\}$.
\end{exm}

\begin{exm}
Again suppose that $n=4$, $\mathrm{dim}Y=5$, and $2m=4$, but now
suppose that $\sigma|_Y$ looks locally like
$$dz_1\wedge dz_2+dz_3\wedge (fdz_4+gdz_5)$$
where $f$ and $g$ are functions on $Y$ whose zero loci meet
transversally. This time $\sigma|_Y$ drops rank along the codimension
two subset $\Sigma:=\{f=0,g=0\}$. Then $F$ is generated by
$g\frac{\partial\phantom{z}}{\partial
  z_4}-f\frac{\partial\phantom{z}}{\partial z_5}$ and hence is locally
free on $Y\backslash\Sigma$, and
$$\omega^{\wedge 2}=dz_1\wedge dz_2\wedge dz_3\wedge (fdz_4+gdz_5)$$
will be a non-trivial section of $\Lambda^4(TY/F)^{\vee}$, and even
non-vanishing on $Y\backslash\Sigma$.
\end{exm}

We now return to the proof of Lemma~\ref{constantrank}.

\vspace*{3mm}
\noindent
{\bf Claim:} A generic $(-1)$-curve $C\cong\P^1$, as found in the previous subsection, 
will not intersect $\Sigma$. In other words, both $F|_C$ and $TY/F|_C$ will be locally
free.

\vspace*{3mm}
Recall that
$\Sigma\subset Y$ has codimension at least two, i.e., dimension at
most $n-1$. The image $\pi(\Sigma)$ under the projection
$\pi:Y\rightarrow\P^n$ will also have dimension at most
$n-1$. Moreover, the subset
$$\pi(\Sigma)_0:=\{t\in\pi(\Sigma)|\Sigma\mbox{ contains the entire curve
}Y_t\}\subset\P^n$$
will have dimension at most $n-2$. It follows that a generic line
$\ell\subset\P^n$ will not intersect $\pi(\Sigma)_0$, and therefore if
$t\in\ell\cap\pi(\Sigma)$ then $Y_t\cap\Sigma\subset Z\cap\Sigma$ will
consist of finitely many points $\{p_1,\ldots,p_k\}$. We have to show that a $(-1)$-curve in $Z$ will not pass through any of these points.

Since $\ell$ is generic, the curve $Y_t$ can have at worst a single node. If one of the points $p_i$ is this node, then a $(-1)$-curve in $Z$ will not pass through $p_i$ by Lemma~\ref{notcontained2}. So we may assume that every point $p_i$ lies in the smooth locus of $Y_t$, and therefore has an associated hyperplane $H_{p_i}\subset\P^n$. By Lemma~\ref{notcontained}, if $\ell\subset\P^n$ does not lie in any of the hyperplanes $H_{p_i}$ then none of the points $p_1,\ldots,p_k$ will be contained in a $(-1)$-curve in $Z$. We must show that this condition can be
satisfied simultaneously for every $t\in\ell\cap\pi(\Sigma)$; for this we use a dimension-counting argument.

For each $t\in\pi(\Sigma)\backslash\pi(\Sigma)_0$, the lines $\ell$ through $t$ contained in one of the hyperplanes $H_{p_i}$ form a finite union of $(n-2)$-dimensional families. Moreover, $\pi(\Sigma)$ itself has
dimension at most $n-1$, so taking the union of these families of
lines as $t$ varies in $\pi(\Sigma)\backslash\pi(\Sigma)_0$ will give a family $\mathcal{S}$
of lines of dimension at most
$$(n-2)+(n-1)=2n-3.$$
The family of all lines in $\P^n$ has dimension $2n-2$, and therefore
a generic line $\ell\subset\P^n$ will not belong to
$\mathcal{S}$. This completes the proof of the claim.

\vspace*{3mm}
Let $C$ be a generic $(-1)$-curve in $Y$; restricting everything to
$C\cong\P^1$ we will show that $2m=2$. In the exact sequence
$$0\rightarrow F|_C\rightarrow TY|_C\rightarrow
TY/F|_C\rightarrow 0$$
the quotient $TY/F|_C$ is locally free of rank $2m$ and therefore
isomorphic to
$$\mathcal{O}(a_1)\oplus\ldots\oplus\mathcal{O}(a_{2m})$$
for some integers $a_1\leq\ldots\leq a_{2m}$. Now
we use the fact that $\omega^{\wedge m}$ is a non-trivial section of
$\Lambda^{2m}(TY/F)^{\vee}$. An argument similar to the one
above shows that this section is still non-trivial when restricted to
$C$, i.e.,
$$\Lambda^{2m} (TY/F)^{\vee}|_C=\mathcal{O}(-a_1-\ldots
-a_{2m})$$
admits a non-trivial section and
we must have $a_1+\ldots +a_{2m}\leq 0$.

In the short exact sequence
$$0\rightarrow TC\rightarrow TY|_C\rightarrow N_{C\subset
  Y}\rightarrow 0$$
the first term $TC$ is isomorphic to $\mathcal{O}(2)$. To identify
  $N_{C\subset Y}$ we use the short exact sequence
$$0\rightarrow N_{C\subset Z}\rightarrow N_{C\subset
  Y}\rightarrow N_{Z\subset Y}|_C\rightarrow 0.$$
Note that $N_{C\subset Z}\cong\mathcal{O}(-1)$ as $C$ is a
$(-1)$-curve in $Z$, and
$$N_{Z\subset
  Y}|_C=\pi^*N_{\ell\subset\P^n}=\pi^*(\mathcal{O}_{\ell}(1)^{\oplus
  (n-1)})=\mathcal{O}(k)^{\oplus (n-1)}$$
where $k\geq 1$ is the degree of the finite covering $C\rightarrow\ell$
induced by $\pi:Z\rightarrow\ell$ (we will see shortly that $k$ must
equal one). Combining the three short exact sequences above gives
$$\begin{array}{ccccccccc}
 & & & & 0 & & & & \\
 & & & & \downarrow & & & & \\
 & & & & \mathcal{O}(2) & & & & \\
 & & & & \downarrow & & & & \\
0 & \rightarrow & F|_C & \rightarrow & TY|_C & \rightarrow &
 \mathcal{O}(a_1)\oplus\ldots\oplus\mathcal{O}(a_{2m}) & \rightarrow &
 0
 \\
 & & & & \downarrow & & & & \\
0 & \rightarrow & \mathcal{O}(-1) & \rightarrow & N_{C\subset Y} &
 \rightarrow & \mathcal{O}(k)^{\oplus (n-1)} & \rightarrow & 0 \\
 & & & & \downarrow & & & & \\
 & & & & 0 & & & & \\
\end{array}$$
One can easily show that $a_1\leq -2$ is impossible; thus there are
two possibilities.

\vspace*{3mm}
\noindent
{\bf Case 1:} Suppose $a_1=-1$. Then $a_2\leq 1$. This is only
possible if $k=1$ and all of the short exact sequences above split, i.e.,
$$N_{C\subset Y}\cong\mathcal{O}(-1)\oplus\mathcal{O}(1)^{\oplus (n-1)}$$
and
$$TY|_C\cong\mathcal{O}(-1)\oplus\mathcal{O}(1)^{\oplus
  (n-1)}\oplus\mathcal{O}(2).$$
Moreover, we must have $2m=2$, and $a_2$ must equal $1$.

\vspace*{3mm}
\noindent
{\bf Case 2:} Suppose $a_1=0$. Then $a_2$ must also equal $0$. This is
only possible if $k=1$,
$$N_{C\subset
  Y}\cong\mathcal{O}^{\oplus 2}\oplus\mathcal{O}(1)^{\oplus
  (n-2)},$$
and
$$TY|_C\cong\mathcal{O}^{\oplus 2}\oplus\mathcal{O}(1)^{\oplus
  (n-2)}\oplus\mathcal{O}(2).$$
Again we must have $2m=2$.

\vspace*{3mm}
In both cases $2m=2$. Thus $\sigma|_Y$ has rank two everywhere and $F$
is locally free of rank $n-1$ everywhere. This completes the proof.
\end{prf}

\begin{cor}
\label{inleaves}
A generic $(-1)$-curve $C$ is contained in a leaf of the foliation
$F$.
\end{cor}

\begin{prf}
In both cases one and two in the proof of Lemma~\ref{constantrank}, we
find that
$$F|_C\cong\mathcal{O}(1)^{\oplus (n-2)}\oplus\mathcal{O}(2).$$
Moreover, the map $TC\rightarrow TY|_C$ lifts to a map $TC\rightarrow F|_C$ which is an isomorphism onto the direct summand $\mathcal{O}(2)$, meaning that $C$ is contained in a leaf of $F$.
\end{prf}

\begin{lem}
\label{simplyconnected}
All leaves of $F$ are rationally connected and therefore simply
connected.
\end{lem}

\begin{prf}
By Lemma~\ref{constantrank}, $F$ is locally free of rank $n-1$, so the
leaves of $F$ are smooth of dimension $n-1$. Since
$$F|_C\cong\mathcal{O}(1)^{\oplus (n-2)}\oplus\mathcal{O}(2)$$
is an ample vector bundle, a generic $(-1)$-curve $C$ satisfies the
criterion
of Theorem~\ref{kebekus} of Kebekus et al. Since $F$ is regular
everywhere, we conclude from the theorem that all leaves of $F$ are
rationally connected. It is well known that smooth rationally
connected varieties are simply connected (see Corollary 4.18 of
Debarre~\cite{debarre01}, for example).
\end{prf}

\subsection{The space of leaves}

A foliation with compact, simply connected leaves will have a
particularly well-behaved space of leaves.

\begin{lem}
\label{hausdorff}
The space of leaves $Y/F$ of the characteristic foliation on $Y$ is a smooth, compact
surface $S$.
\end{lem}

\begin{prf}
Holmann~\cite{holmann80} proved that a holomorphic foliation on a
K{\"a}hler manifold, with all leaves compact, is stable. This means
that every open neighbourhood of a leaf $L$ contains a saturated
neighbourhood of $L$, i.e., a neighbourhood consisting of a union of
leaves. Stability of the foliation is equivalent to the space of
leaves $Y/F$ being Hausdorff.

Let us say a few words about the local structure of the space of
leaves $Y/F$. Let $L$ be a (compact) leaf of $F$, represented by a
point in $Y/F$. Take a small slice $V$ in $Y$ transverse to the
foliation, with $L$ intersecting $V$ at a point $0\in V$. The holonomy
map is a group homomorphism from the fundamental group of $L$ to the
group of automorphisms of $V$ which fix $0$. The holonomy group $H(L)$
of $L$ is the image of the holonomy map. Then $V/H(L)$ is a local model
for the space of leaves $Y/F$ in a neighbourhood of the point
representing $L$ (see Holmann~\cite{holmann78} for details). In our
case, all the leaves are simply connected by
Lemma~\ref{simplyconnected}. Therefore the holonomy groups must be
trivial and all local models for $Y/F$ are smooth.

Since $F$ has rank $n-1$, the space of leaves $Y/F$ is
two-dimensional, and compactness follows from compactness of $Y$.
\end{prf}

\begin{lem}
\label{symplectic}
The space of leaves $S$ admits a holomorphic symplectic form; thus $S$
is either a K3 surface or an abelian surface.
\end{lem}

\begin{prf}
Heuristically, the leaves of $F$ are along the null directions of
$\sigma|_Y$, and therefore $\sigma|_Y$ should descend to $S$ under the
projection $Y\rightarrow S$. To make this rigorous, we will define a
holomorphic symplectic form on $S$ via a local model of $S$, and then
prove that it is independent of the local model chosen.

As explained in the previous proof, a local model for the
space of leaves in a neighbourhood of the point representing the leaf
$L$ is given by a small slice $V_1$ in $Y$, with $L$ intersecting
$V_1$ transversally at a point $0$. In a neighbourhood of $0$, $V_1$
will be transverse to the foliation $F$ and therefore the restriction
$\sigma|_{V_1}$ will be a non-degenerate two-form. Suppose $V_2$ is a
second small slice transverse to $L$, though not necessarily
intersecting $L$ at the same point $L\cap V_1$. Shrinking $V_1$ and
$V_2$ if necessary, there is an isomorphism
$\phi:V_1\rightarrow V_2$ given by taking $L^{\prime}\cap V_1$ to
$L^{\prime}\cap V_2$, where $L^{\prime}$ is an arbitrary leaf of
$F$. This map is well-defined because the
leaves are simply connected, and it takes the point $0$ in $V_1$
(i.e., $L\cap V_1$) to $0$ in $V_2$ (i.e., $L\cap V_2$). In this way,
we identify different local models for the space of leaves.
\begin{figure}
\begin{center}
\includegraphics[width=35.0mm]{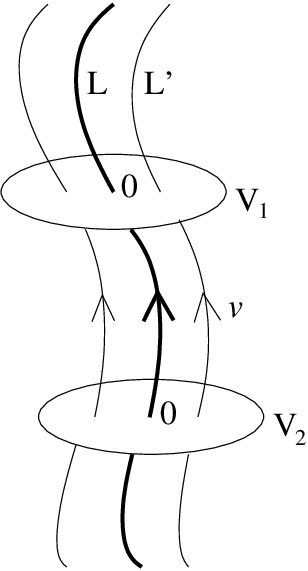}
\end{center}
\caption{Identifying different local models of $Y/F$}
\end{figure}

Observe that $\phi$ can be regarded as the time $t=1$ map of a flow
$\phi_t$ associated to a vector field $v$ along the leaves of
$F$. Since $d\sigma=0$ and $\iota_v(\sigma|_Y)=0$, the Lie derivative
$$\mathcal{L}_v(\sigma|_Y) = \iota_v(d\sigma|_Y)-d(\iota_v\sigma|_Y)$$
vanishes, and therefore the flow $\phi_t$ will preserve the
holomorphic two-form $\sigma|_Y$. In particular
$$\phi^*(\sigma|_{V_2})=\sigma|_{V_1}$$
which completes the proof.
\end{prf}

\begin{lem}
\label{leaves}
Every leaf of $F$ is isomorphic to $\P^{n-1}$. Moreover, the leaves
map isomorphically to hyperplanes in $\P^n$ under the projection
$\pi:Y\rightarrow\P^n$.
\end{lem}

\begin{prf}
Let $C$ be a generic $(-1)$-curve in $Y$, which is contained in a leaf
$L$ by Corollary~\ref{inleaves}. Moreover
$$F|_C\cong\mathcal{O}(1)^{\oplus (n-2)}\oplus\mathcal{O}(2)$$
so that the normal bundle of $C$ inside $L$ is
$$N_{C\subset L}\cong\mathcal{O}(1)^{\oplus (n-2)}.$$
Let $H_L$ be the image of $L$ under the projection $\pi:Y\rightarrow\P^n$. In the proof of Lemma~\ref{constantrank} we saw that $k=1$, i.e., $C\subset L$ projects
isomorphically to $\ell\subset H_L\subset\P^n$ under $\pi$. By Lemma~\ref{transverse} the foliation is transverse to the curves $Y_t$ (at least at non-singular points of $Y_t$), and therefore the morphism
$$d\pi(N_{C\subset L})\cong d\pi(\mathcal{O}_C(1)^{\oplus (n-2)})\rightarrow N_{\ell\subset H_L}\subset N_{\ell\subset\P^n}\cong\mathcal{O}_{\ell}(1)^{\oplus (n-1)}$$
is injective. This forces
$$d\pi(N_{C\subset L})\cong N_{\ell\subset H_L}\cong\mathcal{O}_{\ell}(1)^{\oplus (n-2)},$$
$H_L$ is a hyperplane in $\P^n$, and $\pi$ maps $L$ birationally onto $H_L$.

By Zariski's Main Theorem, the fibres of $\pi:L\rightarrow H_L$ are
connected, as the generic fibre is a single point. Moreover, Lemma~\ref{transverse} implies that the intersection of $L$ with a fibre $Y_t$ of $Y\rightarrow\P^n$ must be zero-dimensional, and therefore $L$ meets every fibre $C_t$ in a single point. This means that $\pi:L\rightarrow H_L$ is a bijection of sets, and thus an isomorphism of varieties since $H_L\cong\P^{n-1}$. This proves the lemma for the leaf $L$ containing $C$,
but by stability of the foliation the result must hold for all leaves.
\end{prf}

\begin{cor}
\label{linearsystem}
Let $Y_t$ be an arbitrary curve in the family $\pi:Y\rightarrow\P^n$. Then $Y_t$
maps isomorphically to its image under the projection $Y\rightarrow
S=Y/F$ to the space of leaves.
\end{cor}

\begin{prf}
The fibres of $Y\rightarrow S=Y/F$ are the leaves of the foliation. By
Lemma~\ref{transverse} the smooth locus of $Y_t$ is everywhere transverse to the leaves,
and by Lemma~\ref{leaves} each leaf will intersect $Y_t$ in at most
one point. The result now follows for smooth curves $Y_t$. If $Y_t$ is singular, we obtain a birational map from $Y_t$ to its image in $S$, which is a bijection on the level of sets. Since $Y_t$ and its image both have (arithmetic) genus $n$, they must in fact be isomorphic.
\end{prf}

We can now complete the proof of Theorem~\ref{main}.

\begin{prf}
By Corollary~\ref{linearsystem}, each curve $Y_t$ maps isomorphically to its image in $S$,
so that $S$ contains an $n$-dimensional linear system of genus $n$
curves. Thus $S$ cannot be an abelian surface, and must be a K3
surface. Therefore $Y\rightarrow\P^n$ is a complete linear
system of genus $n$ curves on a K3 surface $S$. This shows that
$X=\overline{J}^d(Y/\P^n)$ is a Beauville-Mukai
integrable system.
\end{prf}

\section{O'Grady's ten-dimensional space}

In~\cite{ogrady99} O'Grady constructed a new ten-dimensional
holomorphic symplectic manifold, which is not deformation equivalent
to the Hilbert scheme $\mathrm{Hilb}^5S$ of five points on a K3
surface. O'Grady's space may be deformed to a Lagrangian fibration
whose generic fibre is the Jacobian of a smooth genus five curve (see
Rapagnetta~\cite{rapagnetta07}). In this section we describe some features of this Lagrangian fibration, and how it relates to Theorem~\ref{main}.

O'Grady's space $\widetilde{M}^{ss}(2,0,-2)$ is given by desingularizing
the moduli space $M^{ss}(2,0,-2)$ of rank two semi-stable sheaves with
Chern classes $c_1=0$ and $c_2=4$ on a K3 surface $S$. After a
deformation, we may
assume that $S$ is a double cover of the plane branched over a generic
sextic; such a K3 surface is polarized by the pull-back $H$ of a line in the 
plane, which satisfies $H^2=2$. Let $Y\rightarrow |2H|$ be the family of curves 
in the five-dimensional linear system $|2H|$. Curves in this family are pull-backs 
of conics in the plane; the generic curve is a smooth genus five curve, but the family 
also contains reducible curves (pull-backs of pairs of lines) and non-reduced curves
(pull-backs of double lines). Nevertheless, we can {\em define\/} the compactified relative
Jacobian $\overline{J}^d(Y/|2H|)$ as the irreducible component
of the Mukai moduli space of semi-stable sheaves on $S$ which contains
$\iota_*L$, where $L$ is a degree $d$ line bundle on a generic curve $Y_t$ and
$\iota:Y_t\hookrightarrow S$ is the inclusion. In other words
$$\overline{J}^d(Y/|2H|):=M^{ss}(0,2H,d+1-g)=M^{ss}(0,2H,d-4).$$

If $d$ is odd then the Mukai vector $(0,2H,d-4)$ is primitive (note
that $H$ is indivisible for generic $S$); then all semi-stable sheaves
are stable and $M^{ss}(0,2H,d-4)$ is smooth and a deformation of
$\mathrm{Hilb}^5S$. If $d$ is even then $M^{ss}(0,2H,d-4)$ will be
singular along the locus of strictly semi-stable sheaves. In this
case, Rapagnetta showed that $M^{ss}(0,2H,d-4)$ admits a symplectic
desingularization, by using the same construction as for O'Grady's
desingularization of $M^{ss}(2,0,-2)$. Moreover, the desingularization
$\widetilde{M}^{ss}(0,2H,d-4)$ is a Lagrangian fibration over
$|2H|\cong\P^5$.

When $d=6$ let $\iota_*L$ be a generic element of
$\overline{J}^6(Y/|2H|)=M^{ss}(0,2H,2)$, with $L$
a degree six line bundle on a generic (smooth) curve $Y_t$. 
Assume
that $h^0(L)=2$ and that $L$ is globally generated; this behaviour is
generic. Then there is an exact sequence
$$0\rightarrow\mathcal{E}\rightarrow\mathrm{H}^0(L)\otimes\mathcal{O}_S\rightarrow\iota_*L\rightarrow
0$$
and $\mathcal{F}:=\mathcal{E}\otimes\mathcal{O}(H)$ is a stable rank
two bundle in $M^{ss}(2,0,-2)$. In this way, O'Grady described a
birational map between $\overline{J}^6(Y/|2H|)$
and $M^{ss}(2,0,-2)$. Thus the desingularizations
$\widetilde{M}^{ss}(0,2H,2)$ and $\widetilde{M}^{ss}(2,0,-2)$ are birational
and therefore deformation equivalent (see
Huybrechts~\cite{huybrechts97}). In summary, we see that O'Grady's
space $\widetilde{M}^{ss}(2,0,-2)$ can be deformed to the Lagrangian
fibration which comes from desingularizing
$\overline{J}^6(Y/|2H|)$; the
generic fibre is the Jacobian $\mathrm{Pic}^6(Y_t)$ of a smooth genus five curve $Y_t$.

The Lagrangian fibration $\widetilde{M}^{ss}(0,2H,2)$ fails to satisfy a hypothesis of Theorem~\ref{main}; namely, the curves in the family $Y\rightarrow |2H|$ are not all reduced and irreducible. On the other hand, the conclusion still holds: the family of curves {\em is\/} a complete linear system of curves on a K3 surface. Theorem~\ref{main} could possibly be strengthened to allow for reducible and/or non-reduced curves. One expects that the family of curves will still be a linear system of curves on a K3 surface, but the Lagrangian fibration may be a singular Beauville-Mukai system like $\overline{J}^6(Y/|2H|)$, which may or may not admit a symplectic desingularization.

\section{Appendix}

In Lemma~\ref{notminimal} we showed that the surface $Z$ arising in
our construction was not minimal. The hypothesis
$\mathrm{deg}\Delta >4n+20$ was used there to show that $K_Z^2<0$,
implying that $Z$ has Kodaira dimension $-\infty$ if it is minimal. Without this hypothesis, $Z$ could be minimal of Kodaira dimension zero, one, or two. A minimal surface of Kodaira dimension zero with $q=0$ and $p_g=1$ must be a K3 surface; but a K3 surface cannot admit a basepoint free pencil of genus $n\geq 2$ curves. We are left with the following question:

\vspace*{3mm}
\noindent
{\em Given $n\geq 2$, does there exists a minimal surface of Kodaira dimension one or two, with $q=0$ and $p_g=1$, which is fibred over $\P^1$ by genus $n$ curves?\/}

\vspace*{3mm}
\noindent
A negative answer to this question would enable us to remove the hypothesis $\mathrm{deg}\Delta >4n+20$ from Theorem~\ref{main}. We do not know of any such surfaces for $n\geq 4$. However, there do exist such surfaces of general type for $n=2$ and $3$.

\subsection{Kodaira dimension one}

\begin{lem}
If $Z$ is a minimal surface of Kodaira dimension one with $q=0$ and
$p_g=1$, then $Z$ is obtained from an elliptic K3 surface $\tilde{Z}$
by applying logarithmic transforms to smooth fibres and/or fibres of
type $I_k$ ($k\geq 2$).
\end{lem}

\begin{prf}
If $\mathrm{kod}(Z)=1$, then $Z$ must be an elliptic surface over some
curve $B$. A one-form on $B$ would pull back to a one-form on $Z$, but
there are no one-forms on $Z$ as $q=0$, and therefore $B$ is a
rational curve. Suppose that $Z\rightarrow B$ has multiple fibres
$E_i$ with multiplicities $m_i$ (these are not necessarily smooth: they
could be of type $I_k$ with $k\geq2$). Using also $p_g=1$, Kodaira's
formula for the canonical bundle of $Z$ gives
$$K_Z=\sum(m_i-1)E_i.$$
Applying inverse logarithmic transforms to remove the multiple
fibres will produce an elliptic surface $\tilde{Z}$ over $B\cong\P^1$
with trivial canonical bundle. Thus $\tilde{Z}$ must be a K3 surface,
and the lemma follows.
\end{prf}

Unfortunately a logarithmic transform is an analytic, {\em not\/}
algebraic, construction, and so in general we cannot relate curves on
$Z$ to curves on the K3 surface $\tilde{Z}$.

\subsection{Kodaira dimension two}

In this subsection we describe some examples of minimal surfaces of
Kodaira dimension two (general type) with $q=0$ and $p_g=1$, which are
fibred over $\P^1$ by genus $n$ curves. However, we only know of
examples for $n=2$ and $3$.

A bidouble cover is a finite flat morphism with Galois group
$(\mathbb{Z}/2)^2$. In~\cite{catanese99} Catanese constructed various
general type surfaces as bidouble covers of $\P^1\times\P^1$ and
$\P^2$. In particular, consider a smooth bidouble cover $Z$ of
$\P^1\times\P^1$ with branch curves $D_1$, $D_2$, and $D_3$. Then $Z$
has $q=0$ and $p_g=1$ only if the branch curves have bidegree
$$(3,1),\qquad (1,3),\qquad\mbox{and}\qquad (1,1)\qquad\mbox{respectively,}$$
or
$$(4,1),\qquad (0,3),\qquad\mbox{and}\qquad (2,1)\qquad\mbox{respectively,}$$
or
$$(4,0),\qquad (0,4),\qquad\mbox{and}\qquad (2,2)\qquad\mbox{respectively}$$
(these three examples appear on page 103 of
Catanese~\cite{catanese99}). The two projections of $\P^1\times\P^1$
to $\P^1$ then yield fibrations on $Z$. Let us calculate the genus of
the fibres.

In the first example, the fibres will be $4$-to-$1$ covers of $\P^1$
with Galois group $(\mathbb{Z}/2)^2$, and the degrees of the branch
loci will be $3$, $1$, and $1$. We can break this up into two branched
double covers $C_2\rightarrow C_1\rightarrow\P^1$. Then
$C_1\rightarrow\P^1$ will have $1+3=4$ branch points, and by the
Riemann-Hurwitz formula $C_1$ is an elliptic curve. Pulling back the
remaining branch loci, we find that $C_2\rightarrow C_1$ has two
branch points, and therefore $C_2$ has genus two. Therefore the fibres
of $Z\rightarrow\P^1$ have genus two, regardless of which projection
$\P^1\times\P^1\rightarrow\P^1$ we use (the bidegrees are symmetric).

In the second example, a similar calculation shows that the fibres of
$Z\rightarrow\P^1$ have genus two when we project $\P^1\times\P^1$ to
the first factor, and genus three when we project to the second
factor. In the third example the fibres of $Z\rightarrow\P^1$ have
genus three for both projections $\P^1\times\P^1\rightarrow\P^1$.

\begin{flushleft}
Department of Mathematics\hfill sawon@email.unc.edu\\
University of North Carolina\hfill www.unc.edu/$\sim$sawon\\
Chapel Hill NC 27599-3250\\
USA\\
\end{flushleft}

\end{document}